\documentclass[12pt]{amsart}
\usepackage{amsfonts, amsbsy, amsmath, amssymb}

\usepackage[T1]{fontenc}
\usepackage[utf8]{inputenc}
\usepackage{lmodern}
\usepackage{amsmath, amsthm, amssymb,amscd, mathrsfs, amsfonts, mathtools,tikz-cd}
\usepackage{relsize}
\usepackage{euler}
\usepackage{mathtools}

\usepackage{pbox}
\usepackage{times}
\usepackage[all]{xy}
\usepackage{todonotes}
\usepackage{xcolor}
\usepackage[bookmarks=true,
bookmarksnumbered=true, breaklinks=true,
pdfstartview=FitH, hyperfigures=false,
plainpages=false, naturalnames=true,
colorlinks=true,pagebackref=true,
pdfpagelabels]{hyperref}

\usepackage{tikz}
\usepackage{pgfplots}

\pgfplotsset{compat=1.10}
\usepgfplotslibrary{fillbetween}
\usetikzlibrary{patterns}
\usetikzlibrary{matrix}

\usepackage[lite,alphabetic]{amsrefs}

\renewcommand{\PrintDOI}[1]{\href{http://dx.doi.org/\detokenize{#1}}{doi: \detokenize{#1}}%
	\IfEmptyBibField{pages}{, (to appear in print)}{}}

\def\commutatif{\ar@{}[rd]|{\circlearrowleft}}

\newtheorem{theorem}{Theorem}[section]
\newtheorem{prop}[theorem]{Proposition}
\newtheorem{lem}[theorem]{Lemma}
\newtheorem{cor}[theorem]{Corollary}

\theoremstyle{definition}
\newtheorem{definition}[theorem]{Definition}

\theoremstyle{remark}
\newtheorem{rmk}[theorem]{Remark}

\newtheorem{example}[theorem]{Example}

\allowdisplaybreaks

\newcommand{\Z}{\mathbb{Z}}

\usepackage[draft, commentmarkup=todo,
  todonotes={textsize=tiny, textwidth=0.83in}]{changes}
\definechangesauthor[name={Boris}, color=blue]{BT}

\usepackage{bbm}

\def \Z{\mathbb{Z}}

\def \O{\mathcal{O}}

\def \GL{\textup{GL}}

\def \O{\mathcal{O}}

\hypersetup{
	colorlinks = true,
	urlcolor = blue,
	linkcolor = blue,
	citecolor = red,
	pdfpagemode = UseNone
}

\input{t1ptm.fd}
\DeclareFontShape{T1}{ptm}{m}{up}{<->ssub*ptm/m/n}{}
\begin{document}

\title{A Helly theorem for the Bruhat-Tits building}

\begin{abstract}
Let $K$ be a discretely valued field with valuation ring $\O$. We show that a finite family of $\mathcal{O}$-lattices in $K^m$ is adapted to a common basis if any subfamily of at most $m+1$ $\mathcal{O}$-lattices is. 
In other words, in the Bruhat-Tits building of $\GL(m)$, a finite family of simplices lies in a common apartment if and only if any subfamily of at most $m+1$ simplices does. The proof uses a version of Radon's theorem for $\mathcal{O}$-lattices. We also give an alternative proof using truncated modules.
\end{abstract}

\author{Rudrendra Kashyap}
\address{Department of Mathematics, University of Pittsburgh, Pittsburgh, PA, USA.}
\email{ruk26@pitt.edu}

\author{Kiumars Kaveh}
\address{Department of Mathematics, University of Pittsburgh, Pittsburgh, PA, USA.}
\email{kaveh@pitt.edu}

\author{Christopher Manon}
\address{Department of Mathematics, University of Kentucky, Lexington, KY, USA}
\email{Christopher.Manon@uky.edu}

\author{Boris Tsvelikhovskiy} 
\address{Department of Mathematics, University of California, Riverside, CA, USA} 
\email{borist@ucr.edu}

\maketitle

\section{Introduction}

Let $W$ be a finite-dimensional vector space over a field $\Bbbk$. A collection of flags $\mathcal{F}_1,  \dots, \mathcal{F}_n$ in $W$ is said to be \textit{adapted} to a basis $\mathcal{B}$ of $W$ if each subspace appearing in any of the flags $\mathcal{F}_i$ is spanned by a subset of $\mathcal{B}$. A remarkable theorem of Klyachko \cite[Theorem~6.1.2]{Klyachko} establishes a local-to-global criterion for the simultaneous adaptation of flags:

\begin{theorem}
\label{thm:klyachko-helly-flags}
Let $S = \{\mathcal{F}_1, \dots, \mathcal{F}_n\}$ be a collection of flags in an $m$-dimensional vector space $W$. For $S$ to be simultaneously adapted to a single basis $\mathcal{B}$ of $W$, it is necessary and sufficient that every subcollection of at most $m+1$ flags from $S$ be adapted to some basis of $W$.
\end{theorem}

Following Klyachko, we refer to Theorem~\ref{thm:klyachko-helly-flags} as a \textit{Helly-type theorem} for vector subspaces and flags. It serves as a direct linear algebraic analogue to the classical Helly theorem in convex geometry, which asserts that a finite family of convex sets in $\mathbb{R}^m$ has non-empty intersection if and only if every subcollection of at most $m+1$ sets has non-empty intersection.

In this paper, we extend Klyachko's result to lattices over discrete valuation rings and give its 
building-theoretic formulation.

\subsection*{Modules over Discrete Valuation Rings}
Let $\mathcal{O}$ be a discrete valuation ring (DVR) with fraction field $K$, and let $\varpi$ be a uniformizer for $\mathcal{O}$ (i.e., a generator of its unique maximal ideal $\mathfrak{m}$). 

Let $E$ be an $m$-dimensional $K$-vector space. By an $\mathcal O$-lattice in $E$ we mean a finitely generated
$\mathcal O$-submodule $\Lambda\subseteq E$ such that $K\Lambda=E$. Since $\mathcal O$ is a DVR, every such module is free of rank $m$.  Throughout we may refer to an $\mathcal{O}$-lattice simply as a lattice. A family of $\mathcal{O}$-lattices $\{\Lambda_i\}_{i \in I}$ in $E$ is said to be \textit{(simultaneously) adapted} to a $K$-basis $\mathcal{B} = \{e_1, \dots, e_m\}$ of $E$ if each lattice $\Lambda_i$ decomposes as
\[
  \Lambda_i = \bigoplus_{j=1}^m \varpi^{a_{ij}} \mathcal{O} e_j
\]
for some $a_{ij} \in \mathbb{Z}$.

We establish an exact lattice analogue of Klyachko's theorem:

\begin{theorem}[See Theorem~\ref{thm-main}]
\label{thm:lattice-helly-intro}
Let $S = \{\Lambda_1, \Lambda_2, \dots, \Lambda_n\}$ with $n \ge m+1$ be a collection of $\mathcal{O}$-lattices of rank $m$ in an $m$-dimensional $K$-vector space $E$. If any $m+1$ lattices from $S$ are simultaneously adapted to some basis of $E$, then there exists a basis of $E$ to which all lattices in $S$ are simultaneously adapted.
\end{theorem}

A key step in the proof of this result is an $\mathcal{O}$-lattice version of the tropical Radon theorem in \cite[Section~3]{GM}. Namely, for $m+1$ lattices already adapted to one basis, we can find two disjoint nonempty subfamilies of these $m+1$ lattices whose suitably scaled sums coincide (Theorem \ref{thm:radon_lattices}).

The lattice version of Radon's theorem enables us to choose a
reference lattice $\Lambda_0$. After picking a homothetic copy
$L$ of $\Lambda_0$ that contains all the given lattices, we associate to each lattice an increasing filtration of $L/\varpi L$. Klyachko's theorem gives a basis of this residue space adapted to all the filtrations. A lifting argument then produces an $\mathcal O$-basis of $L$ to which all the original lattices are adapted.

We point out that the threshold $m+1$ is optimal for $m\geq2$: a minimally incompatible configuration of $m+1$ lines in the residue space gives a family of $m+1$ lattices such that every proper subfamily admits a common adapted basis, but the entire family does not.

In Section~\ref{sec:alternative}, we also give an alternative proof of the lattice Helly theorem (Theorem~\ref{thm-main}). Unlike the proof given in Section~\ref{sec:lattices}, which uses a direct lifting argument and a lattice version of Radon's theorem, the alternative proof establishes the conditional reference-lattice theorem by induction on $N$ for submodules of free $\mathcal O/(\varpi^N)$-modules, using intersections of affine subspaces to construct compatible lifts. 

For $K=k((t))$ and $\mathcal O=k[[t]]$, the lattice theorem has a direct interpretation on the affine Grassmannian of $\operatorname{GL}_m$, whose $k$-points parametrize full $\mathcal O$-lattices in $K^m$. A basis apartment consists of all lattices adapted to a fixed $K$-basis. Thus a finite collection of points lies in one basis apartment whenever every subcollection of at most $m+1$ points does, and this threshold is optimal for $m\geq2$. The same implication holds for affine partial flag varieties, whose points are periodic lattice chains: if every subcollection of at most $m+1$ chains
admits a basis adapted to all its constituent lattices, then one basis adapts all the chains. This follows by applying the lattice theorem to the finite collection of lattices occurring in one period of each chain.

In Section~\ref{sec:buildings}, we reinterpret Theorem \ref{thm:lattice-helly-intro} as a combinatorial property of the Bruhat-Tits building in type A. In the type A Bruhat-Tits building, vertices are homothety classes of lattices. A collection of vertices lies in a common apartment if and only if the corresponding lattices are simultaneously adapted to a basis. 



\subsection*{Related Work}
As mentioned earlier, this work is motivated by \cite{Klyachko}. Huang and the second author \cite[\S6]{HK} introduced Helly numbers of buildings, motivated by the splitting problem for toric principal bundles, and asked for sharp bounds. In a recent paper \cite{Tsvelikhovskiy2026}, the last author obtains exact values and bounds for the Helly numbers of connected reductive complex algebraic groups and applies these results to the equivariant splitting of toric
principal bundles. Finally, in \cite{JK}, Jafari and the second author extend Klyachko's theorem to the flats of matroids.

\subsection*{Organization of the Paper}
In Section~\ref{sec:review} we recall Klyachko's theorem and its formulation in terms of parabolic subgroups. Section~\ref{sec:lattices} contains the proof of the lattice theorem using a reference lattice, reduction to the residue field, and a lattice analogue of the tropical Radon theorem. 
Section~\ref{sec:buildings} interprets the results in the language of buildings. Section~\ref{sec:alternative} gives
a second proof using truncated modules and minimum-weight assignments.

\subsection*{Acknowledgment}
While this paper is written by the authors, AI tools were utilized to assist in improving the exposition and reviewing the literature. In addition, the statement and proof of Theorem~\ref{thm:radon_lattices} were proposed by OpenAI ChatGPT (GPT-5.6 Sol). The material was polished, reviewed and verified by the authors, who assume full responsibility for the manuscript.


\section{The \texorpdfstring{$\operatorname{GL}(V)$}{GL(V)} Case: Review}
\label{sec:review}

Let $\Bbbk$ be a field and $V$ an $m$-dimensional vector space over $\Bbbk$. A collection of subspaces $\{W_1, W_2, \dots, W_n\}$ of $V$ is said to be \emph{adapted} to a basis $\mathcal{B}$ of $V$ if each subspace $W_i$ is spanned by a subset of $\mathcal{B}$. A remarkable result of Klyachko \cite[Theorem~6.1.2]{Klyachko} gives a local-to-global criterion for having a common adapted basis:

\begin{theorem}
\label{thm:Klyachko}
Let $\mathcal{S} = \{W_1, \dots, W_n\}$ with $n\ge m+1$ be a family of subspaces in an $m$-dimensional vector space $V$. The collection $\mathcal{S}$ is adapted to a common basis of $V$ if and only if every subcollection of $m+1$ subspaces from $\mathcal{S}$ is adapted to some basis of $V$.
\end{theorem}

Theorem \ref{thm:Klyachko} can be extended to flags of subspaces in $V$. A flag of subspaces in $V$ is adapted to a basis $\mathcal{B}$ if all the subspaces in it are adapted to $\mathcal{B}$. Given a finite family of flags in $V$, applying Theorem \ref{thm:Klyachko} to the family of all the subspaces appearing in these flags, we obtain the following: 
A finite family of flags in $V$ is adapted to a common basis if and only if any  subfamily of at most $m+1$ flags is. 

The above statement in turn can be formulated in terms of parabolic subgroups in  the general linear group $\operatorname{GL}(V)$ as follows. A choice of basis $\mathcal{B}$ for $V$ uniquely determines a maximal torus $T \subset \operatorname{GL}(V)$ consisting of transformations that are diagonal with respect to the basis $\mathcal{B}$. A flag $\mathcal{F}$ in $V$ is adapted to $\mathcal{B}$ if and only if its stabilizer $P = \operatorname{Stab}_{\operatorname{GL}(V)}(\mathcal{F})$ contains $T$. Since parabolic subgroups are exactly stabilizers of flags, we have the following:

\begin{theorem}
\label{thm:parabolics}
Let \(V\) be an \(m\)-dimensional vector space over an algebraically closed field \(\Bbbk\), with \(m\geq 1\), and let \(\{P_\alpha\}_{\alpha\in A}\) be a finite family of parabolic subgroups of \(\operatorname{GL}(V)\). Then \(\bigcap\limits_{\alpha\in A}P_\alpha\) contains a maximal torus if and only if the intersection of every subfamily of at most \(m+1\) members contains a maximal torus.
\end{theorem}

\section{A Helly Theorem for Free Modules over a DVR}
\label{sec:lattices}

 In this section, we extend Theorem~\ref{thm:Klyachko} to the setting of free modules over a discrete valuation ring. Let $\mathcal{O}$ be a discrete valuation ring with fraction field $K$, and let $\varpi$ be a uniformizer for $\mathcal{O}$ (i.e., a generator of its unique maximal ideal $\mathfrak{m}$). 
Write $\Bbbk=\mathcal O/(\varpi)$ for the residue field.

\begin{definition}
    Let \(E = K^m\) be an \(m\)-dimensional vector space over \(K\). An \textbf{\(\mathcal{O}\)-lattice} \(\Lambda\) in \(E\) is a rank \(m\) (free) \(\mathcal{O}\)-submodule of \(E\). We say that a lattice \( \Lambda \subset E \) is \emph{adapted} to a basis \( f_1, \dots, f_m \) of \( E \) if \( \Lambda \) has an \(\mathcal{O} \)-basis of the form \( \varpi^{a_1} f_1, \ldots, \varpi^{a_m} f_m \) for some integers \( a_1, \dots, a_m \in \mathbb{Z}\).
\end{definition}

\begin{theorem}
\label{thm-main}
Let  \(S = \{\Lambda_1, \Lambda_2, \dots, \Lambda_n\}\) with $n\ge m+1$ be a collection of \(\mathcal{O}\)-lattices in \( K^m \). If every subfamily of at most \(m+1\) lattices from this collection is adapted to some basis of \( K^m \), then there exists a basis of \(K^m \) to which all lattices in the collection are simultaneously adapted.
\end{theorem}

The proof of Theorem~\ref{thm-main} proceeds in two steps. We first establish a conditional variant, Theorem~\ref{th-main-body}, which proves the result under the assumption that there exists a fixed reference lattice $\Lambda_0$ such that $\Lambda_0$ and every $(m+1)$-element subset of $S$ are simultaneously adapted to a common basis. We then demonstrate that the auxiliary lattice $\Lambda_0$ can be chosen from the geometry of the given collection $S$ without imposing restrictive hypotheses on $S$. 


\subsection{Step 1: Conditional Result Assuming a Reference Lattice}

We start with a technical lemma.

\begin{lem}\label{lem:finite-intersections}
Let $\Lambda_1,\ldots,\Lambda_s$ be lattices in an $m$-dimensional $K$-vector space. There are subsets $I,J\subseteq\{1,\ldots,s\}$, each of cardinality at most $m$, such that
\[
  \sum_{i=1}^s\Lambda_i=\sum_{i\in I}\Lambda_i,
  \qquad  \bigcap\limits_{i=1}^s\Lambda_i=\bigcap\limits_{i\in J}\Lambda_i.
\]
Moreover, if $M$ is a finite-length $\mathcal O$-module with
\[
  \dim_\Bbbk\ker(\varpi:M\to M)=m,
\]
every finite intersection of submodules of $M$ is the intersection of at most $m$ of them.
\end{lem}

\begin{proof}
Put $\Lambda:=\sum\limits_i\Lambda_i$. The images of the elements of $\bigcup\limits_i\Lambda_i$ span $\Lambda/\varpi \Lambda$. Select $m$ such elements whose images form a basis. By Nakayama's lemma, they form an $\mathcal O$-basis of $\Lambda$. They belong to at most $m$ of the $\Lambda_i$, proving the assertion about sums.

For a lattice $A\subseteq E$, write
\[
  A^\vee=\{f\in E^*:f(A)\subseteq\mathcal O\}.
\]
Apply the sum assertion to the $\Lambda_i^\vee$ and take duals, using
\[
  A^{\vee\vee}=A,\qquad
  \bigcap\limits_{i=1}^s\Lambda_i
  =\left(\sum_{i=1}^s\Lambda_i^\vee\right)^\vee
  =\left(\sum_{i\in J}\Lambda_i^\vee\right)^\vee
  =\bigcap\limits_{i\in J}\Lambda_i.
\]

For the last assertion, the structure theorem over a DVR gives a surjection $p:\mathcal O^m\to M$. Since $\mathcal O^m/\ker p\simeq M$ has finite length, $\ker p$ is a lattice in $K^m$. The inverse images of the submodules are therefore lattices. Apply the lattice assertion and pass to the quotient by $\ker p$.
\end{proof}

\begin{theorem} \label{th-main-body}
    Let $\Lambda_1, \Lambda_2, \dots, \Lambda_n$ with $n\ge m+1$ be a collection of $\mathcal{O}$-lattices of rank $m$ in the vector space $E=K^m$. 
    Suppose there exists a rank $m$ $\mathcal{O}$-lattice $\Lambda_0$ in $E$ such that for every subset $\mathcal{Y} \subset \{\Lambda_1, \Lambda_2, \dots, \Lambda_n\}$ of size $|\mathcal{Y}| = m+1$, there is a basis $\mathcal{B}_{\mathcal{Y}}$ of $E$ with the property that each lattice in $\mathcal{Y}\cup\{\Lambda_0\}$ is adapted to $\mathcal{B}_{\mathcal{Y}}$.  
    Then there exists a basis $\mathcal{B}$ of $E$ to which all lattices $\Lambda_0,\Lambda_1, \Lambda_2, \dots, \Lambda_n$ are simultaneously adapted.
\end{theorem}

\begin{proof}
Choose a homothetic copy $L = \varpi^{-N}\Lambda_0$ containing all $\Lambda_i$ for $1 \le i \le n$, and set $$\overline L = L/\varpi L.$$ For $q \in \Z_{\geq 0}$, define
\begin{equation}
    \label{eq:filtration_of_L}
    F_i(q) = \frac{(L \cap \varpi^{-q}\Lambda_i) + \varpi L}{\varpi L}.
\end{equation}
These form increasing filtrations of $\Bbbk$-subspaces of $\overline L$, which stabilize at $\overline L$ for sufficiently large $q$.

By hypothesis, for every subfamily of at most $m+1$ lattices from $\{\Lambda_1, \dots, \Lambda_n\}$, there exists a common adapted basis with $\Lambda_0$ (and thus with $L$). This induces a basis of $\overline L$ splitting the corresponding filtrations. 

The collection of subspaces appearing across the filtrations \eqref{eq:filtration_of_L} is finite. Since any subcollection of at most $m+1$ such subspaces arises from at most $m+1$ distinct filtrations, they admit a common adapted basis. Consequently, Theorem~\ref{thm:Klyachko} applied over the residue field $\Bbbk$ yields a basis $\{\bar{f}_1, \dots, \bar{f}_m\}$ of $\bar{L}$ that is simultaneously adapted to all filtrations.

For each $1 \le j \le m$, define
\[
    a_{ij} \coloneqq \min\{q \geq 0 \mid \overline f_j \in F_i(q)\} \quad \text{and} \quad A_{ij} \coloneqq L \cap \varpi^{-a_{ij}}\Lambda_i.
\]
By Lemma~\ref{lem:finite-intersections}, there exists a subset $J_j \subset \{1, \dots, n\}$ with $|J_j| \leq m$ such that
\[
    \bigcap\limits_{i=1}^n A_{ij} = \bigcap\limits_{i\in J_j} A_{ij}.
\]
As $|J_j| \leq m$, the reference-lattice hypothesis supplies a common adapted basis for $L$ and the $\Lambda_i$ with $i\in J_j$. In this basis the $A_{ij}$, for $i\in J_j$, are coordinate
sublattices of $L$. Therefore the reduction map
$p:L\to\bar L$ satisfies
\begin{equation}
\label{eq:proj_intersection_commute}
  p\left(\bigcap\limits_{i=1}^n A_{ij}\right)
  =p\left(\bigcap\limits_{i\in J_j}A_{ij}\right)
  =\bigcap\limits_{i\in J_j}p(A_{ij}).
\end{equation}

By definition of $a_{ij}$,
$\bar f_j\in F_i(a_{ij})=p(A_{ij})$ for every $i$.
Choose
\[
  f_j\in\bigcap\limits_{i=1}^n A_{ij},
  \qquad p(f_j)=\bar f_j.
\]
The vectors $f_1,\ldots,f_m$ form an $\mathcal O$-basis of
$L$ by Nakayama's lemma. Moreover,
$f_j\in L\cap\varpi^{-a_{ij}}\Lambda_i$ implies
$\varpi^{a_{ij}}f_j\in\Lambda_i$. Thus
\[
  L_i'=\bigoplus_{j=1}^m\varpi^{a_{ij}}\mathcal O f_j
  \subseteq\Lambda_i.
\]
Let $b_{i1}, \dots, b_{im} \ge 0$ be the elementary-divisor exponents of $\Lambda_i$ in $L$. The dimension of $F_i(q)$ satisfies
\[
    \dim_\Bbbk F_i(q) = \#\{j \mid b_{ij} \leq q\} = \#\{j \mid a_{ij} \leq q\}.
\]
Thus the multisets $\{a_{i1}, \dots, a_{im}\}$ and $\{b_{i1}, \dots, b_{im}\}$ coincide. Therefore,
\[
    \operatorname{length}_{\mathcal O}(L/L_i') = \sum_{j=1}^m a_{ij} = \sum_{j=1}^m b_{ij} = \operatorname{length}_{\mathcal O}(L/\Lambda_i).
\]
Since $L_i' \subseteq \Lambda_i \subseteq L$, exactness of length implies $\operatorname{length}_{\mathcal O}(\Lambda_i/L_i') = 0$, so $L_i' = \Lambda_i$. Hence $\{f_1, \dots, f_m\}$ is a simultaneously adapted basis for all $\Lambda_1, \dots, \Lambda_n$.

Furthermore, the basis adapts $\Lambda_0$, as it forms an $\mathcal{O}$-basis for $L = \varpi^{-N}\Lambda_0$.
\end{proof}

\subsection{Step 2: The Choice of Auxiliary Reference Lattice (Radon's Theorem for Lattices)}
\label{subsec:Radon_for_lattices}

Recall Radon's classical theorem, which asserts that any set of $m+2$ points in $\mathbb{R}^m$ can be partitioned into two disjoint subsets whose convex hulls intersect. Below, we give an $\mathcal{O}$-lattice analogue over a discrete valuation ring $\mathcal{O}$ with uniformizer $\varpi$. We use an integral lattice formulation of the max-plus Radon theorem \cite[Section~3]{GM}. Relative to a common adapted basis, lattice sums become coordinatewise minima of exponent vectors, and homotheties become translations by multiples of $(1,\ldots,1)$. This result serves as a key step in the proof of Theorem~\ref{thm-main} and is also of independent interest.

\begin{theorem}[Lattice Radon's Theorem]
\label{thm:radon_lattices}
Let $m \ge 1$, let $E$ be an $m$-dimensional $K$-vector space, and let $\Lambda_1, \dots, \Lambda_{m+1} \subset E$ be lattices adapted to a common basis. Then there exist disjoint, non-empty index sets $I,J$ with 
\[
\{1,\ldots,m+1\}=I\sqcup J,
\qquad I\neq\varnothing,\quad J\neq\varnothing
\]
and integer shifts $\{c_k\}_{k \in I \cup J} \subset \mathbb{Z}$ such that
\begin{equation}
\label{eq:radon_lattice_sum}
    \sum_{i \in I} \varpi^{c_i} \Lambda_i = \sum_{j \in J} \varpi^{c_j} \Lambda_j.
\end{equation}
\end{theorem}

\begin{proof}
Since the lattices $\Lambda_1, \ldots, \Lambda_{m+1}$ are adapted to the basis $\{e_1, \ldots, e_m\}$, each lattice decomposes as:
\[
    \Lambda_i = \bigoplus_{r=1}^m \varpi^{a_{ri}} \mathcal{O} e_r, \qquad \text{for } i = 1, \ldots, m+1,
\]
where $a_{ri} \in \mathbb{Z}$ represents the exponent of $\varpi$ along the basis vector $e_r$.

To construct the index partition and integer scaling factors, consider the formal matrix 
\[
A(t) = (t^{a_{ri}}) \in \operatorname{Mat}_{m \times (m+1)}(\mathbb{Q}(t))
\]
over the  field $\mathbb{Q}(t)$ of rational functions in one variable. As $A(t)$ has size $m \times (m+1)$, its kernel is non-trivial. Choose a non-zero vector 
\[
u(t) = (u_1(t), \ldots, u_{m+1}(t))^T \in \ker A(t).
\]

For each index $i$ in the non-zero support $\operatorname{supp}(u)$, let $c_i = \operatorname{ord}_{t=0} u_i(t) \in \mathbb{Z}$ be the order of vanishing at $t=0$, and write 
\[
u_i(t) = \kappa_i t^{c_i} + O(t^{c_i+1})
\]
with leading coefficient $\kappa_i \in \mathbb{Q} \setminus \{0\}$. Partition the support into two disjoint sets based on the signs of the leading coefficients:
\[
    I \coloneqq \{ i \mid \kappa_i > 0 \}, \qquad J \coloneqq \{ i \mid \kappa_i < 0 \}.
\]
Since $u(t) \in \ker A(t)$, for each row $r \in \{1, \ldots, m\}$ we have the identity
\[
    \sum_{i \in I \cup J} t^{a_{ri}} u_i(t) = 0 \in \mathbb{Q}(t).
\]
For the lowest power of $t$ in this sum to have a zero coefficient, the minimum valuation attained among the positive terms ($i \in I$) must match the minimum valuation attained among the negative terms ($j \in J$). Consequently, for every coordinate $r \in \{1, \ldots, m\}$, we have:
\begin{equation}
\label{eq:coordinate_exponent_equality}
    \min_{i \in I} (a_{ri} + c_i) = \min_{j \in J} (a_{rj} + c_j).
\end{equation}
Both $I$ and $J$ are necessarily non-empty to satisfy this cancellation across all rows.

Finally, observing that the $r$-th coordinate of the lattice sum $\sum\limits_{i \in I} \varpi^{c_i} \Lambda_i$ is generated by $\varpi^{\min\limits_{i \in I} (a_{ri} + c_i)} \mathcal{O} e_r$, the coordinate-wise exponent equalities in \eqref{eq:coordinate_exponent_equality} yield the explicit module identity in \eqref{eq:radon_lattice_sum}.

Let $\Lambda_*$ denote the common lattice sum from \eqref{eq:radon_lattice_sum}. For each unused index $k\in\{1,\ldots,m+1\}\setminus(I\cup J)$, choose an
integer $c_k$ large enough that $\varpi^{c_k}\Lambda_k
\subseteq\Lambda_*$, and add $k$ to $I$. The equality is unchanged, and $I,J$ now partition all indices.
\end{proof}

\begin{rmk}
\label{rmk:radon-number}
The threshold $m+1$ is directly consistent with the classical Radon number from convex geometry. 

Indeed, relative to a fixed adapted basis of $K^m$, an  $\mathcal{O}$-lattice $\Lambda$ is uniquely determined by its tuple of valuation exponents
\[
    (a_1, \dots, a_m) \in \mathbb{Z}^m.
\]
Scaling the lattice by $\varpi^c$ shifts its exponent vector to $(a_1 + c, \dots, a_m + c) = (a_1, \dots, a_m) + c(1, \dots, 1)$. Consequently, homothety classes of adapted lattices correspond to points in the discrete quotient $\mathbb{Z}^m / \mathbb{Z}(1, \dots, 1)$, which embeds naturally into the real quotient space
\[
    \mathbb{R}^m / \mathbb{R}(1, \dots, 1).
\]
Since $\dim_{\mathbb{R}} \left(\mathbb{R}^m / \mathbb{R}(1, \dots, 1)\right) = m - 1$, the ambient geometric space has dimension $d = m - 1$. Since the classical Radon number of a $d$-dimensional affine space is $d + 2$, the corresponding threshold in this lattice setting is precisely
\[
    (m - 1) + 2 = m + 1.
\]
\end{rmk}

\subsection{Completing the Proof of Theorem \ref{thm-main}}

We now combine the results of the two preceding subsections to establish the main assertion and thereby complete the proof of Theorem \ref{thm-main}.

\begin{proof}[Proof of Theorem \ref{thm-main}]
The case $m=0$ is immediate. We thus assume $m\geq1$. 
We proceed by induction on the size $n$ of the collection $S$. The case $n=m+1$ is exactly the hypothesis. Assume $n\geq m+2$ and that the assertion holds for families of $n-1$ lattices.

By hypothesis, $\Lambda_1,\ldots,\Lambda_{m+1}$ admit a
common adapted basis. Apply Theorem~\ref{thm:radon_lattices} to these lattices
and write
\[
  \Lambda_0:=\sum_{i\in I}\varpi^{c_i}\Lambda_i
           =\sum_{j\in J}\varpi^{c_j}\Lambda_j,
  \qquad I\cap J=\varnothing,
  \quad I,J\neq\varnothing.
\]

For each $k \in \{1, \dots, n\}$, consider the subset $\mathcal{Y}_k := S \setminus \{\Lambda_k\}$ of size $n-1$. By the inductive hypothesis, all lattices in $\mathcal{Y}_k$ are simultaneously adapted to a common basis $\mathcal{B}_k$. 

Since $I \cap J = \varnothing$, the index $k$ cannot belong to both $I$ and $J$, so at least one index set corresponds to lattices contained entirely in $\mathcal{Y}_k$. Expressing $\Lambda_0$ through that sum represents it as an $\mathcal{O}$-linear sum of homothetic lattices in $\mathcal{Y}_k$, guaranteeing that $\Lambda_0$ is also adapted to $\mathcal{B}_k$.

Now, let $\mathcal{Y}' \subset S$ be an arbitrary subset of size $|\mathcal{Y}'| = m+1$. Because $n \ge m+2$, we have $m+1 \le n-1$, so $\mathcal{Y}' \subseteq \mathcal{Y}_k$ for at least one $k \in \{1, \dots, n\}$. Setting $\mathcal{B}_{\mathcal{Y}'} := \mathcal{B}_k$, it follows that every lattice in $\mathcal{Y}' \cup \{\Lambda_0\} \subseteq \mathcal{Y}_k \cup \{\Lambda_0\}$ is adapted to $\mathcal{B}_{\mathcal{Y}'}$.

Since $\mathcal{Y}'$ was an arbitrary subset of size $m+1$, all hypotheses of Theorem \ref{th-main-body} are satisfied for the collection $S$ and the lattice $\Lambda_0$. Applying Theorem \ref{th-main-body} yields a basis of $K^m$ to which all lattices $\Lambda_1, \dots, \Lambda_n \in S$ are simultaneously adapted, completing the inductive step.
\end{proof}

 The following example shows that the threshold $m+1$ cannot be replaced by $m$ when $m\geq2$.

\begin{example}
\label{ex:lattice-helly-sharpness}
Let $m \ge 2$, let $\{e_1, \dots, e_m\}$ be a $K$-basis of $K^m$, and consider the reference lattice
\[
  \Lambda_0 = \bigoplus_{i=1}^{m} \mathcal{O} e_i.
\]
Fix a vector $v = \sum\limits_{i=1}^{m} \alpha_i e_i$ with unit coefficients $\alpha_i \in \mathcal{O}^\times$, and define a collection of $m+1$ lattices by
\[
  \Lambda_i = \mathcal{O} e_i + \varpi \Lambda_0 \quad (1 \le i \le m),
  \qquad
  \Lambda_{m+1} = \mathcal{O} v + \varpi \Lambda_0.
\]

We claim that every subcollection of $m$ lattices from $\{\Lambda_1, \dots, \Lambda_{m+1}\}$ is adapted to a common basis, but the full collection of $m+1$ lattices is not.

To see this for any subcollection of size $m$:
\begin{itemize}
  \item If $\Lambda_{m+1}$ is omitted, the remaining $m$ lattices $\{\Lambda_1, \dots, \Lambda_m\}$ are manifestly adapted to the standard basis $\{e_1, \dots, e_m\}$.
  \item If $\Lambda_j$ is omitted for some $1 \le j \le m$, consider the basis 
  \[
  \mathcal{B}_j = \{e_1, \dots, e_{j-1}, v, e_{j+1}, \dots, e_m\}.
  \]
  As $\alpha_j \in \mathcal{O}^\times$, $\mathcal{B}_j$ forms an $\mathcal{O}$-basis of $\Lambda_0$, and all $m$ lattices in $\{\Lambda_1, \dots, \Lambda_{m+1}\} \setminus \{\Lambda_j\}$ are adapted to $\mathcal{B}_j$.
\end{itemize}

Now suppose, toward a contradiction, that all $m+1$ lattices were simultaneously adapted to some basis $\mathcal{B}'$ of $K^m$. Then their sum
\[
  \sum_{i=1}^{m} \Lambda_i = \Lambda_0
\]
would also be adapted to $\mathcal{B}'$. By rescaling the basis elements of $\mathcal{B}'$ by suitable powers of $\varpi$, we may assume without loss of generality that $\mathcal{B}'$ is an $\mathcal{O}$-basis of $\Lambda_0$. 

Reducing modulo $\varpi \Lambda_0$, the image of each lattice $\Lambda_i$ ($1 \le i \le m+1$) in the $m$-dimensional residue vector space $\bar{\Lambda}_0 = \Lambda_0 / \varpi \Lambda_0 \cong \Bbbk^m$ (where $\Bbbk = \mathcal{O}/\varpi\mathcal{O}$) is a 1-dimensional subspace:
\[
  \bar{\Lambda}_i = \Bbbk \bar{e}_i \quad (1 \le i \le m),
  \qquad
  \bar{\Lambda}_{m+1} = \Bbbk \bar{v}.
\]
Since $m\ge2$ and every $\alpha_i$ is a unit, the lines
$\Bbbk\bar e_1,\ldots,\Bbbk\bar e_m,\Bbbk\bar v$ are pairwise distinct. An $m$-dimensional vector space has only $m$ coordinate lines for a fixed basis, so these $m+1$ lines cannot all be coordinate, yielding a contradiction. 

Consequently, the threshold $m+1$ in Theorem~\ref{thm-main} is sharp.
\end{example}

\begin{rmk}
Theorem~\ref{thm-main} is an affine-building
analogue of Klyachko's Helly theorem \cite[Theorem~6.1.2]{Klyachko}, rather than a formal consequence of that theorem over the fraction field $K$ (see Section \ref{sec:buildings}). 
From the point of view of algebraic groups, Theorem \ref{thm-main} concerns the (affine) Bruhat-Tits building and parahoric subgroups, the stabilizers of $\mathcal O$-lattices in $\operatorname{GL}_m(K)$, whereas Klyachko's theorem concerns the (spherical) Tits building and  parabolic subgroups, the stabilizers of flags in $\operatorname{GL}_m(K)$.

We also remark that the rank-two case of Theorem \ref{thm-main} was previously noted in \cite[Example~5.4]{KMT}. 
\end{rmk}

\section{Interpretation in Terms of Buildings}
\label{sec:buildings}

We recall the building-theoretic language needed to interpret the
preceding results; see \cite{AB} for general background. A
building is a simplicial complex equipped with subcomplexes called
\emph{apartments}, each isomorphic to a fixed Coxeter complex, such
that any two simplices lie in an apartment and any two apartments
admit an isomorphism fixing their intersection pointwise.
Its maximal simplices are called \emph{chambers}. There are two important classes of buildings: \emph{spherical buildings} and \emph{affine buildings}. 
Spherical apartments are finite Coxeter complexes; affine apartments
are infinite Euclidean Coxeter complexes.

\begin{definition}
\label{defn:apartment-helly}
Let $\Delta$ be a building equipped with a specified apartment system. The \emph{Helly number} $h(\Delta)$ is the minimal positive integer $h$ such that any finite family of simplices $\mathcal{S}$ is contained in a single apartment provided every subfamily of at most $h$ simplices is contained in an apartment. If no such integer exists, we set $h(\Delta) = \infty$. 
\end{definition}

The above definition goes back to \cite[Section~6]{HK}.
For a simplex $\sigma$, we define the set of apartments containing $\sigma$ by
\[
  \mathcal{A}(\sigma) \coloneqq \{ A \mid A \text{ is an apartment and } \sigma \subseteq A \}.
\]
Then a family $\mathcal S$ lies in one apartment precisely when
\[
  \bigcap\limits_{\sigma\in\mathcal S}\mathcal A(\sigma)\neq\varnothing.
\]
Thus $h(\Delta)$ is the ordinary Helly number of this incidence family in the set of apartments.

\subsection{Spherical buildings and common maximal tori}

Let $G$ be a connected reductive complex algebraic group of positive semisimple rank. The spherical building $\Delta_{\mathrm{sph}}(G)$ has simplices corresponding to proper parabolic subgroups of $G$, ordered by reverse inclusion. Under this identification, maximal proper parabolic subgroups correspond to vertices, and Borel subgroups correspond to chambers. Each maximal torus $T \subseteq G$ determines the apartment
\[
  A_T \coloneqq \{ P \subsetneq G \mid P \text{ is parabolic and } T \subseteq P \}.
\]
Consequently, parabolic subgroups $P_1,\ldots,P_N$ correspond to simplices in one apartment if and only if the intersection $\bigcap\limits_iP_i$ contains a maximal torus of $G$. An occurrence of $P_i=G$ imposes no condition and may
be omitted.

For $G=\operatorname{GL}(V)$, Theorem~\ref{thm:parabolics} 
shows that the Helly number of the spherical building of $G$ is equal to $\dim V+1$.

\subsection{Affine buildings and lattices}

Let $\mathcal O$ be a discrete valuation ring with fraction field $K$ and uniformizer $\varpi$, and let $E$ have dimension $m\geq2$ over $K$. We consider the affine building (also called the \emph{Bruhat-Tits building}) $\Delta_{\mathrm{aff}}(E)$ of $\operatorname{PGL}(E)$, equivalently the reduced affine building of $\operatorname{GL}(E)$, with its standard
apartment system. Its vertices are homothety classes $[\Lambda]$
of $\mathcal O$-lattices in $E$.

A collection of $q+1$ vertices form a $q$-simplex if they admit representatives $\Lambda_0, \ldots, \Lambda_q$ satisfying
\[
  \varpi\Lambda_0\subsetneq\Lambda_q\subsetneq\cdots
  \subsetneq\Lambda_1\subsetneq\Lambda_0.
\]
Chambers have $m$ vertices. Repeating such a chain
under multiplication by powers of $\varpi$ gives a \emph{periodic lattice chain}.

A \emph{frame} is an unordered decomposition
\[
  E=\ell_1\oplus\cdots\oplus\ell_m,
  \qquad \dim_K\ell_j=1.
\]
Choosing $e_j\in\ell_j\setminus\{0\}$, the corresponding apartment has vertices
\[
  \left[
    \bigoplus_{j=1}^m\varpi^{a_j}\mathcal Oe_j
  \right],
  \qquad (a_1,\ldots,a_m)\in\mathbb Z^m.
\]
It is a full subcomplex of the building. Its vertices are
parametrized by
\[
  \mathbb Z^m/\mathbb Z(1,\ldots,1),
\]
and its geometric realization is
$\mathbb R^m/\mathbb R(1,\ldots,1)$, with the affine Coxeter triangulation of type $\widetilde A_{m-1}$.
Rescaling the $e_j$ translates the coordinates; it does not change the apartment. 


Our main theorem (Theorem \ref{thm-main}) gives the Helly number of the affine building $\Delta_{\mathrm{aff}}(E)$: 
\begin{theorem}[Helly Numbers of Affine Buildings]
\label{thm:affine-apartment-helly}
For $m \ge 2$, the Helly number of the affine building $\Delta_{\mathrm{aff}}(E)$ is:
\[
  h\bigl(\Delta_{\mathrm{aff}}(E)\bigr) = m + 1.
\]
\end{theorem}

\begin{proof}
Let $\mathcal{S}$ be a finite family of simplices in $\Delta_{\mathrm{aff}}(E)$ such that every subfamily of at most $m+1$ simplices lies in a common apartment. For a simplex $\sigma$, write $V(\sigma)$ for its set of vertices, and set
\[
  V(\mathcal S)=\bigcup_{\sigma\in\mathcal S}V(\sigma).
\]
to be the set of all vertices appearing in these simplices. Any subcollection of at most $m+1$ vertices in $V(\mathcal{S})$ can be chosen from at most $m+1$ simplices in $\mathcal{S}$, and therefore lies in a common apartment by hypothesis. Applying Theorem \ref{thm-main} to $V(\mathcal{S})$ yields a single apartment $\mathcal{A}$ containing all vertices in $V(\mathcal{S})$. Since apartments are full subcomplexes of $\Delta_{\mathrm{aff}}(E)$, $\mathcal{A}$ contains every simplex $\sigma \in \mathcal{S}$. In light of Example~\ref{ex:lattice-helly-sharpness}, we get the equality  $h\bigl(\Delta_{\mathrm{aff}}(E)\bigr) = m+1$.

\end{proof}

\begin{rmk}
\label{rmk:residue-obstruction}
The link of a vertex $[L]$ in the affine building $\Delta_{\mathrm{aff}}(E)$ is naturally isomorphic to the spherical building of the residue space $L/\varpi L \cong \Bbbk^m$: a proper non-zero subspace $U \subset L/\varpi L$ corresponds to the lattice $\Lambda_U \subset L$ containing $\varpi L$ and satisfying $\Lambda_U / \varpi L = U$.

Under this correspondence, the lattices in Example~\ref{ex:lattice-helly-sharpness} originate from a collection of $m+1$ lines in an $m$-dimensional residue space $\Bbbk^m$. We have $L=\Lambda_0 = \Lambda_1 + \dots + \Lambda_m$. Consequently, any apartment containing the corresponding $m+1$ vertices would necessarily contain the vertex $[\Lambda_0]$ as well, and would induce a common coordinate frame for the $m+1$ lines in the residue space $\Lambda_0/\varpi \Lambda_0$. 
\end{rmk}

\subsection{Implications for Affine Grassmannians and Partial Flag Varieties}
\label{subsec:aff_gr_and_flags}

Let $\Bbbk$ be a field, set $\mathcal{O} = \Bbbk[[t]]$ with fraction field $K = \Bbbk((t))$, and let $G = \operatorname{GL}_m$. The set of $\Bbbk$-points of the affine Grassmannian $\operatorname{Gr}_G$ is given by the quotient
\[
  \operatorname{Gr}_G(\Bbbk) = G(K) / G(\mathcal{O}),
\]
which is naturally identified with the set of all  $\mathcal{O}$-lattices in $K^m$ via the map $g G(\mathcal{O}) \mapsto g \cdot \mathcal{O}^m$ (see \cite[Section~1.1]{Zhu17}).

Let $\varnothing \neq J \subseteq \mathbb{Z}/m\mathbb{Z}$, and let $\widetilde{J} \subset \mathbb{Z}$ be its preimage under the quotient map $\mathbb{Z} \to \mathbb{Z}/m\mathbb{Z}$. A \textit{periodic $J$-chain} of $\mathcal{O}$-lattices in $K^m$ is a family $\mathcal{L} = (\Lambda_a)_{a \in \widetilde{J}}$ of  $\mathcal{O}$-lattices satisfying
\[
  \Lambda_b \subset \Lambda_a \quad (a < b), \qquad \dim_{\Bbbk}(\Lambda_a / \Lambda_b) = b - a, \qquad \Lambda_{a+m} = t \Lambda_a.
\]
Fix the standard periodic chain defined, for $a = qm + r$ with $0 \le r < m$, by
\[
  \Lambda_a^{\mathrm{std}} = t^q \left( \bigoplus_{j=1}^r t \mathcal{O} e_j \oplus \bigoplus_{j=r+1}^m \mathcal{O} e_j \right),
\]
and define the associated parahoric subgroup
\[
  \mathcal{P}_J = \{ g \in G(K) \mid g \Lambda_a^{\mathrm{std}} = \Lambda_a^{\mathrm{std}} \text{ for all } a \in \widetilde{J} \}.
\]
Then the set of $\Bbbk$-points 
\[
\operatorname{Fl}_J(\Bbbk) = G(K) / \mathcal{P}_J
\]
is naturally identified with the space of periodic $J$-chains in $K^m$.

For a $K$-basis $\mathcal B=\{b_1,\ldots,b_m\}$ of $K^m$, define
\[
  \mathcal A_{\mathcal B,J}
  \coloneqq
  \left\{
    \mathcal L\in\operatorname{Fl}_J(\Bbbk)
    \;\middle|\;
    \text{every lattice of $\mathcal L$ is adapted to $\mathcal B$}
  \right\}.
\]
We call these subsets \emph{basis apartments}.
For each $\mathcal L\in\mathcal A_{\mathcal B,J}$,
the homothety classes of the lattices in one period form
a simplex in the building apartment determined by
$\mathcal B$.

\begin{rmk}
    Unlike the normalized $\operatorname{SL}_m$ setting, the present $\operatorname{GL}_m$ convention imposes no fixed restriction on the valuation of the determinant.
\end{rmk}

In the special case $J = \{0\}$, the affine partial flag variety $\operatorname{Fl}_J(\Bbbk)$ reduces to the affine Grassmannian $\operatorname{Gr}_G(\Bbbk)$.

For a detailed treatment of parahoric quotients and lattice-chain descriptions of affine flag varieties, see \cite[Sections~2.1.6--2.1.7]{Yun17}.

Theorem~\ref{thm-main} yields a direct Helly-type property for points in affine partial flag varieties.

\begin{cor}[Affine Partial Flag Varieties]
\label{cor:affine-flag-helly}
Let $J \subseteq \mathbb{Z}/m\mathbb{Z}$ be a non-empty subset, and let $x_1, \dots, x_n \in \operatorname{Fl}_J(\Bbbk)$ with $n \ge m+1$. If every subcollection of at most $m+1$ points is contained in some basis apartment $\mathcal{A}_{\mathcal{B}, J}$, then all $n$ points are simultaneously contained in a single basis apartment $\mathcal{A}_{\mathcal{B}, J}$.
\end{cor}

\begin{proof}
Take the finite union $S$ of the lattices in one period
of each chain. Any at most $m+1$ members of $S$ occur
in at most $m+1$ of the chains, so they admit a common
adapted basis by hypothesis. Theorem~\ref{thm-main}
assures the existence of a basis adapted to all of the lattices in $S$. Multiplication by
powers of $t$ preserves adaptation, so this basis adapts
every lattice in all the chains.
\end{proof}

\section{Alternative Proof of the Main Theorem}
\label{sec:alternative}

We give an alternative proof of Theorem~\ref{thm-main}.
The proof in Section~\ref{sec:lattices} uses
Lemma~\ref{lem:finite-intersections} to lift a common adapted
basis from the residue space, and constructs
the auxiliary reference lattice using lattice Radon theorem
(Theorem~\ref{thm:radon_lattices}).
In contrast, here we first prove a simultaneous adaptation theorem for submodules of free $\mathcal O/(\varpi^N)$-modules by induction on $N$, starting from Klyachko's theorem at $N=1$ and constructing compatible lifts through intersections of affine subspaces. We then construct the auxiliary reference lattice using minimum weight assignments instead of the lattice Radon theorem, and complete the proof of Theorem~\ref{thm-main} by induction on the number of lattices.


For $N\geq1$, put $R_N=\mathcal O/(\varpi^N)$. 
A submodule $U$ of a free $R_N$-module $M$ is \emph{adapted} to a basis $e_1,\ldots,e_m$ of $M$ if
it is generated by $\varpi^{a_1}e_1,\ldots,\varpi^{a_m}e_m$ for some integers $0\leq a_j\leq N$.
Equivalently,
\[
  U=\bigoplus_{j=1}^m\varpi^{a_j}R_Ne_j.
\] Here $a_j=N$ gives a zero summand, since $\varpi^NR_N=0$.

\begin{prop}\label{prop:alt-finite-module}
Let $M$ be a free $R_N$-module of rank $m$. A finite family of submodules of $M$ admits a common adapted basis if every subfamily of at most $m+1$ members does.
\end{prop}
\begin{proof}
Denote the given family by $\{U_1,\ldots,U_s\}$. If $s\leq m+1$, the family already admits a common
adapted basis by hypothesis. Thus, we may assume
that $s\geq m+2$. We proceed by induction on $N$.
For $N=1$, the assertion follows from Theorem~\ref{thm:Klyachko}. Suppose $N\geq2$ and assume the proposition for $N-1$.

We first establish the identity 
\begin{equation}
\label{eq:alt-image-intersection}
 \varpi\Bigl(\bigcap\limits_i U_i\Bigr)=\bigcap\limits_i\varpi U_i.
\end{equation}
The inclusion from left to right holds for arbitrary submodules. To prove the reverse inclusion, let $f\in\bigcap\limits_i\varpi U_i$, and set
\[
  X_f\coloneqq\{x\in M\mid\varpi x=f\},
  \qquad
  {C_{f,i}}\coloneqq X_f\cap U_i.
\]

Put
\[
  V=\ker(\varpi:M\to M)=\varpi^{N-1}M.
\]
Since $\varpi V=0$, the $R_N$-module structure on $V$ factors through $\Bbbk=R_N/(\varpi)$.
Thus $V$ is a $\Bbbk$-vector space, and
$\dim_{\Bbbk}V=m$ because $M$ is free of rank $m$
over $R_N$.

Choose $x_0\in X_f$. Such an element exists because
$f\in\varpi U_1\subseteq\varpi M$.
Then
\[
  X_f=x_0+V,
\]
so translation by $-x_0$ identifies $X_f$ with the vector space $V$. In particular, $X_f$ is an $m$-dimensional affine space over $\Bbbk$.

For each $i$, choose $x_i\in C_{f,i}$, which is possible because $f\in\varpi U_i$. We have
\[
  C_{f,i}=x_i+(U_i\cap V).
\]
Since $U_i\cap V$ is a $\Bbbk$-linear subspace of $V$, each $C_{f,i}$ is a nonempty affine subspace of $X_f$.

Let $I\subseteq\{1,\ldots,s\}$ be a nonempty subset of cardinality  $|I|\leq m+1$. By hypothesis, there is a basis $g_1,\ldots,g_m$ of $M$ adapted to all $U_i$ with $i\in I$.
Write
\[
  U_i=\bigoplus_{j=1}^m\varpi^{a_{ij}}R_Ng_j
  \quad (i\in I),
  \qquad
  b_j=\max_{i\in I}a_{ij}.
\]

In this common adapted basis, intersections can be computed coordinatewise:
\[
  \bigcap_{i\in I}U_i
  =\bigoplus_{j=1}^m\varpi^{b_j}R_Ng_j.
\]
Similarly,
\[
  \bigcap_{i\in I}\varpi U_i
  =\bigoplus_{j=1}^m\varpi^{b_j+1}R_Ng_j
  =\varpi\Bigl(\bigcap_{i\in I}U_i\Bigr),
\]
where powers of $\varpi$ with exponent at least $N$ give the zero ideal.

Since $f\in\bigcap_{i\in I}\varpi U_i$, the preceding equality gives an element
$x\in\bigcap_{i\in I}U_i$ satisfying $\varpi x=f$. Thus $x\in C_{f,i}$ for every $i\in I$, and hence
\[
  \bigcap_{i\in I}C_{f,i}\neq\varnothing.
\]

Let $\mathcal C\coloneqq\{C_{f,1},\ldots,C_{f,s}\}$. Among all nonempty intersections of subfamilies of $\mathcal C$, choose one of minimal dimension and denote it by $D$. We regard the intersection of the empty subfamily as $X_f$. Choose a subfamily of $\mathcal C$ whose intersection is $D$.

Starting from $X_f$, intersect successively with the
members of this chosen subfamily, omitting any member that contains the current intersection. Every intermediate intersection contains $D$ and is therefore nonempty.

Each subspace added during this process strictly decreases the dimension.  Since $\dim X_f=m$, we obtain
\[
  D=\bigcap_{i\in J}C_{f,i}
  \qquad\text{for some }J\subseteq\{1,\ldots,s\}
  \text{ with }|J|\leq m.
\]
For every $i$, the intersection $D\cap C_{f,i}$ is nonempty, because it involves at most $m+1$ members of the family.
By the minimality of $\dim D$, it has the same
dimension as $D$. Since it is an affine subspace of $D$, this implies
\[
  D\cap C_{f,i}=D.
\]
Consequently,
\[
  \bigcap_{i=1}^s C_{f,i}\supseteq D\neq\varnothing.
\]

Any element of this intersection belongs to every $U_i$ and maps to $f$ under multiplication by $\varpi$. This proves \eqref{eq:alt-image-intersection}.

The same argument applies to every nonempty subfamily of $\{U_1,\ldots,U_s\}$.

The module $\varpi M$ is free of rank $m$ over $R_{N-1}$.
If $U_i=\bigoplus_j\varpi^{a_{ij}}R_Ng_j$ in an adapted
basis of $M$, then $\varpi U_i$ and $U_i\cap\varpi M$
are adapted to the basis $\varpi g_1,\ldots,\varpi g_m$
of $\varpi M$, with exponents $\min(a_{ij},N-1)$ and
$\max(a_{ij}-1,0)$, respectively.
Every subfamily of at most $m+1$ members of
\[
  \{\varpi U_i,\ U_i\cap\varpi M:1\leq i\leq s\}
\]
involves at most $m+1$ original indices and therefore admits
a common adapted basis. The inductive hypothesis gives an
$R_{N-1}$-basis $f_1,\ldots,f_m$ of $\varpi M$ such that
\[
  \varpi U_i=\bigoplus_j\varpi^{c_{ij}}R_{N-1}f_j,
  \qquad
  U_i\cap\varpi M=\bigoplus_j\varpi^{d_{ij}}R_{N-1}f_j,
  \qquad 0\leq c_{ij},d_{ij}\leq N-1.
\]
We claim that
\[
  (\varpi U_i)\cap\varpi^2M
  =\varpi(U_i\cap\varpi M).
\]
The inclusion from right to left is immediate.
Conversely, write an element of the left-hand side as
$z=\varpi u=\varpi^2v$, where $u\in U_i$ and $v\in M$.
Then
\[
  u-\varpi v\in\ker(\varpi:M\to M)
  =\varpi^{N-1}M\subseteq\varpi M,
\]
where the last inclusion uses $N\geq2$.
Thus $u\in U_i\cap\varpi M$, proving the claim.

Together with
\[
  \varpi(U_i\cap\varpi M)
  \subseteq\varpi U_i\subseteq U_i\cap\varpi M,
\]
this identity shows that, for each $(i,j)$, exactly one of
the following holds:

\begin{equation}\label{eq:alt-exponent-cases}
 (c_{ij},d_{ij})=(0,0),\qquad c_{ij}=d_{ij}+1,
 \qquad (c_{ij},d_{ij})=(N-1,N-1).
\end{equation}
Indeed, the inclusions give $d_{ij}\leq c_{ij}\leq d_{ij}+1$.
If $c_{ij}=d_{ij}=t$ with $0<t<N-1$, then
$\varpi^tf_j$ belongs to $(\varpi U_i)\cap\varpi^2M$
but not to $\varpi(U_i\cap\varpi M)$, contradicting the
identity above.

For each $j$, apply \eqref{eq:alt-image-intersection} to
the subfamily of those $U_i$ for which $c_{ij}=d_{ij}=0$.
Since $f_j\in\varpi U_i$ for every such $i$, there is a
lift $e_j\in M$ satisfying $\varpi e_j=f_j$ and belonging
to every member of this subfamily.
If the subfamily is empty, choose any lift of $f_j$.
Multiplication by $\varpi$ induces an isomorphism
\[
  M/\varpi M\xrightarrow{\sim}\varpi M/\varpi^2M.
\]
The classes of the $e_j$ therefore form a basis of
$M/\varpi M$, so $e_1,\ldots,e_m$ form a basis of $M$
by Nakayama's lemma.

Define $a_{ij}$ to be $0$, $c_{ij}$, or $N$ in the three
respective cases of \eqref{eq:alt-exponent-cases}, and set
\[
  D_i=\bigoplus_j\varpi^{a_{ij}}R_Ne_j.
\]
Then $D_i\subseteq U_i$.
In the first case, this follows from the choice of $e_j$.
In the second case,
\[
  \varpi^{c_{ij}}e_j
  =\varpi^{d_{ij}}f_j\in U_i\cap\varpi M.
\]
In the third case, the generator is zero. In all three cases,
\[
  c_{ij}=\min(a_{ij},N-1),
  \qquad
  d_{ij}=\max(a_{ij}-1,0).
\]
The coordinate formulas give
\[
  \varpi D_i=\varpi U_i,
  \qquad D_i\cap\varpi M=U_i\cap\varpi M.
\]
For $u\in U_i$, choose $d\in D_i$ with $\varpi d=\varpi u$.
Then
\[
  u-d\in U_i\cap\ker(\varpi:M\to M)
  \subseteq U_i\cap\varpi M
  =D_i\cap\varpi M\subseteq D_i.
\]
Hence $u\in D_i$, and therefore $U_i=D_i$ for every $i$.
\end{proof}

\begin{proof}[Alternative proof of Theorem~\ref{th-main-body}]
Choose $c\in\mathbb Z$ such that $L=\varpi^c\Lambda_0$
contains every $\Lambda_i$, and choose $N\geq1$ such that
\[
  \varpi^NL\subseteq\Lambda_i\subseteq L
  \qquad\text{for every }i.
\]
Set $M=L/\varpi^NL$ and $U_i=\Lambda_i/\varpi^NL$. A common adapted basis for $L$ and any chosen subfamily of the $\Lambda_i$ can be rescaled to an $\mathcal O$-basis of $L$. Its images in $M$ form an adapted $R_N$-basis for the
corresponding $U_i$. The hypothesis therefore allows us to apply Proposition~\ref{prop:alt-finite-module} to the family of all $U_i$. Lift the resulting basis to $b_1,\ldots,b_m\in L$. These form an $\mathcal O$-basis of $L$ by Nakayama's lemma. If the exponents of $U_i$ are $a_{ij}\in\{0,\ldots,N\}$, taking inverse images gives $\Lambda_i=\bigoplus_j\varpi^{a_{ij}}\mathcal O b_j$. Since $\Lambda_0=\varpi^{-c}L$, the same basis is
adapted to $\Lambda_0$.
\end{proof}

It remains to choose the reference lattice without increasing the
number of lattices to be checked. The following construction uses the
cofactor identity from tropical linear algebra; compare
\cite[Section~5]{RST}. We include its elementary proof, which does not require a genericity assumption.

\begin{lem}\label{lem:alt-redundant-sum}
If $\Lambda_1,\ldots,\Lambda_{m+1}$ admit a common adapted basis, there
are integers $c_1,\ldots,c_{m+1}$ such that, for every $k$,
\begin{equation}
\label{eq:alt-redundant-sum}
 \sum_{i=1}^{m+1}\varpi^{c_i}\Lambda_i
   =\sum_{\substack{1\leq i\leq m+1\\i\neq k}}
                       \varpi^{c_i}\Lambda_i.
\end{equation}
\end{lem}
\begin{proof}
For $m=0$, both sides are zero, so assume $m\geq1$. Choose a common adapted basis $e_1,\ldots,e_m$ and write
\[
\Lambda_i=\bigoplus_{r=1}^m
\varpi^{a_{ri}}\mathcal O e_r.
\]
For each $i$, define
\[
c_i=\min_{\sigma}\sum_{q=1}^m a_{q,\sigma(q)},
\qquad
\sigma:\{1,\ldots,m\}\xrightarrow{\sim}
\{1,\ldots,m+1\}\setminus\{i\}.
\]

Fix $r$, choose $i$ minimizing $a_{ri}+c_i$, and choose
a bijection $\sigma$ attaining $c_i$. Put $j=\sigma(r)$, so $j\neq i$. Replacing $\sigma(r)=j$ by $i$, while leaving the other values unchanged, gives a bijection onto
$\{1,\ldots,m+1\}\setminus\{j\}$.
Hence
\[
  c_j\leq c_i-a_{rj}+a_{ri},
  \qquad\text{and therefore}\qquad
  a_{rj}+c_j\leq a_{ri}+c_i.
\]
By the choice of $i$, equality holds. Thus $\min_i(a_{ri}+c_i)$ is attained at two distinct indices and is unchanged when any one index is omitted.
Since
\[
  \sum_{i=1}^{m+1}\varpi^{c_i}\Lambda_i
  =\bigoplus_{r=1}^m
  \varpi^{\min_i(a_{ri}+c_i)}\mathcal O e_r,
\]
the same coordinate formula after omitting $k$ proves
\eqref{eq:alt-redundant-sum}.

\end{proof}

\begin{proof}[Alternative proof of Theorem~\ref{thm-main}]

The case $m=0$ is immediate, so assume $m\geq1$.
We proceed by induction on the number $n$ of lattices. If $n\leq m+1$, the assertion follows directly from the hypothesis.

Suppose $n\geq m+2$ and assume the result holds for
families of $n-1$ lattices. By hypothesis, the first $m+1$ lattices admit a common adapted basis. Applying
Lemma~\ref{lem:alt-redundant-sum}, we obtain integers $c_1,\ldots,c_{m+1}$ such that the lattice 
$
\Lambda_0:=\sum_{i=1}^{m+1}\varpi^{c_i}\Lambda_i
$
satisfies 
\[
\Lambda_0
=\sum_{\substack{1\leq i\leq m+1\\i\neq k}}
\varpi^{c_i}\Lambda_i
\qquad (1\leq k\leq m+1).
\]
For each $k\in\{1,\ldots,n\}$, the family
$(\Lambda_i)_{i\neq k}$ satisfies the hypothesis
of the theorem. By induction, there exists a basis
$\mathcal B_k$ adapted to all lattices in this family. 
We claim that $\mathcal B_k$ also adapts $\Lambda_0$.
If $k\leq m+1$, this follows from \eqref{eq:alt-redundant-sum}; if $k>m+1$, it follows from the definition of $\Lambda_0$. Indeed, in both cases $\Lambda_0$ is a sum of homothetic copies of lattices adapted to $\mathcal B_k$, and such sums remain adapted to the same basis.

Let $T\subseteq\{1,\ldots,n\}$ satisfy $|T|\leq m+1$.
Since $n\geq m+2$, choose $k\notin T$.
The basis $\mathcal B_k$ adapts $\Lambda_0$ and every
$\Lambda_i$ with $i\in T$. Thus the hypotheses of
Theorem~\ref{th-main-body} hold. Consequently, there is a common adapted basis
for $\Lambda_0,\Lambda_1,\ldots,\Lambda_n$, completing
the induction.

\end{proof}

\end{document}